\documentclass[10pt,twoside]{cras_i-en}
\usepackage{epsfig}
\usepackage{amssymb,amsbsy,amsmath,amsfonts,amssymb,amscd}
\usepackage[english,francais]{babel}
\usepackage{times}
\input utile.def
\newtheorem{e-proposition}[theorem]{Proposition}

\newtheorem{e-definition}[theorem]{D\'efinition\rm}

\ComParit{S0764-4442}
\PIT{FLA}
\PXMA{????}
\Add{?}
\Volume{XXX}
\Year{XXXX}
\FirstPage{XX}
\LastPage{XX}

\AuteurCourant{J. Casado-D\'{\i}az et al.}
\TitreCourant{An elastic beam fixed on a small part of one of its extremities} 
\Journal
\Rubrique{Probl\`emes math\'ematiques de la m\'ecanique}{Mathematical Problems in Mechanics}
\SousRubrique{}{}
\PresentePar{XXX}{}
\Recu{}{}{}
\begin{document}
\selectlanguage{english}
\title{%
Asymptotic behavior of an elastic beam fixed on a small part of
one of its extremities\\
\medskip
(Preliminary version of a Note to be published in a slightly abbreviated form in
\\
\smallskip
C. R. Acad. Sci. Paris, Ser. I, 338 (2004), pp. 975-980)
}
\author{%
Juan CASADO-D\'IAZ$^{\text{a}}$,\ \ Manuel
LUNA-LAYNEZ$^{\text{a}}$, \ \ Fran\c{c}ois MURAT$^{\text{b}}$}
\address{%
\begin{itemize}\labelsep=2mm\leftskip=-5mm
\item[$^{\text{a}}$]
Departamento de Ecuaciones Diferenciales y An\'alisis Num\'erico,
Universidad de Sevilla, c/ Tarfia s/n, 41012 Sevilla, Espa\~{n}a
\item[$^{\text{b}}$]
Laboratoire Jacques-Louis Lions, Universit\'e Pierre et Marie
Curie, bo\^{i}te courrier 187, 75252 Paris cedex 05, France
\end{itemize}
}
\maketitle

\thispagestyle{empty}
\begin{Abstract}{
We study the asymptotic behavior of the solution of an
anisotropic, heterogeneous, li\-nea\-rized elasticity problem in a
cylinder whose diameter~$\varepsilon$ tends to zero. The cylinder
is assu\-med to be fixed (homogeneous Dirichlet boundary
condition) on the whole of one of its extremities, but only on a
small part (of size $\varepsilon r^\varepsilon$) of the second
one; the Neumann boundary condition is assumed on the remainder of
the boundary. We show that the result depends on~$r^\varepsilon$,
and that there are~3 critical sizes, namely
$r^\varepsilon=\varepsilon^3$, $r^\varepsilon=\varepsilon$,
and~$r^\varepsilon=\varepsilon^{1/3}$, 
and in total~7
different regimes. We also prove a corrector result for each
behavior of~$r^\varepsilon$.}\end{Abstract}
\selectlanguage{french}
\begin{Ftitle}{%
Comportement asymptotique d'une poutre \'elastique fix\'ee sur une
petite partie de l'une de ses extr\'emit\'es}\end{Ftitle}
\begin{Resume}{Nous \'etudions le comportement asymptotique de la
solution d'un probl\`eme d'\'elasticit\'e lin\'eaire anisotrope et
h\'et\'erog\`ene dans un cylindre dont le diam\`etre $\varepsilon$
tend vers z\'ero. Le cylindre est fix\'e (condition de Dirichlet
homog\`ene) sur la totalit\'e de l'une de ses extr\'emit\'es, mais
seulement sur une petite partie (de taille $\varepsilon
r^\varepsilon$) de l'autre base; sur le reste de la fronti\`ere on
a la condition de Neumann. Nous montrons que le r\'esultat depend
de $r^\varepsilon$, 
et qu'il existe 3 tailles critiques, \`a savoir
$r^\varepsilon=\varepsilon^3$, $r^\varepsilon=\varepsilon$ et
$r^\varepsilon=\varepsilon^{1/3}$, 
et au total~7 comportements diff\'erents.
Nous donnons un r\'esultat de
correcteur pour tous les comportements
de~$r^\varepsilon$.}\end{Resume}

\AFv
Dans cette Note, nous \'etudions le comportement asymptotique de
la solution d'un probleme d'\'elasticit\'e lin\'eaire anisotrope
et h\'et\'erog\`ene pos\'e dans un cylindre
$\Omega^\varepsilon=(0,1)\times \varepsilon S\subset \RR^3$ dont
le diam\`etre $\varepsilon$ tend vers z\'ero et dont l'axe est le
premier axe de coordonn\'ees~$Ox_1$. A l'une de ses extr\'emit\'es
($x_1=1$) le cylindre est fix\'e (condition de Dirichlet
homog\`ene) sur la totalit\'e de la base
$\Gamma^\varepsilon_1=\{1\}\times \varepsilon S$, tandis qu'\`a
l'autre extr\'emit\'e ($x_1=0$), il est seulement fix\'e sur une
petite partie $\Gamma^\varepsilon_0=\{0\}\times \varepsilon
r^\varepsilon S_0$ de la base, partie qui est bien plus petite
que~$\varepsilon$ car $r^\varepsilon$ tend vers z\'ero. Sur le
reste de la fronti\`ere de $\Omega^\varepsilon$, on a la condition
de Neumann (voir Figure~1). Pour $A$ tenseur d'ordre~4 coercif \`a
coefficients continus dans $\overline \Omega$ ($\Omega=(0,1)\times
S$), et $f\in L^2(\Omega)^3$, $h\in L^2(\Omega)^{3\times 3}_s$
donn\'es, on d\'efinit $A^\varepsilon$, $F^\varepsilon$ et
$H^\varepsilon$ par~(\ref{defAeps}) et~(\ref{defFepsHeps}) et on
d\'efinit $U^\varepsilon$ comme la solution du probl\`eme
d'\'elasticit\'e~(\ref{probdehom}).
\par Le but de cette Note est de d\'ecrire le comportement
asymptotique de $U^\varepsilon$ et de donner un r\'esultat de
correcteur pour $e(U^\varepsilon)$ quand $\varepsilon$ et
$r^\varepsilon$ tendent vers z\'ero. Ce r\'esultat est d\'ecrit
dans le Th\'eor\`eme enonc\'e dans la version anglaise ci-dessous,
qui fait apparaitre 3 tailles critiques, \`a savoir
$r^\varepsilon\approx\varepsilon^3$,
$r^\varepsilon\approx\varepsilon$ et
$r^\varepsilon\approx\varepsilon^{1/3}$, 
qui s\'eparent~4 zones correspondant \`a 
$r^\varepsilon\ll \varepsilon^3$, $\varepsilon^3\ll
r^\varepsilon \ll\varepsilon$, $\varepsilon\ll r^\varepsilon \ll
\varepsilon^{1/3}$ et $\varepsilon^{1/3}\ll r^\varepsilon\le C$,
donc au total 7 cas diff\'erents. Le r\'esultat de convergence est
donn\'e par~(\ref{correctorUeps}), o\`u $u$ est d\'efini par la
solution de~(\ref{18dicpbmalim}), tandis
que~(\ref{correctore(Ueps)}) est un r\'esultat de correcteur pour
$e(U^\varepsilon)$.
\par Quand $r^\varepsilon\ll\varepsilon^3$, l'ensemble
$\Gamma^\varepsilon_0=\{0\}\times \varepsilon r^\varepsilon S_0$
est si petit que la condition de Dirichlet impos\'ee \`a
$U^\varepsilon$ pour $x_1=0$ dispara\^{\i}t compl\`etement \`a la
limite. Quand $\varepsilon^3\ll r^\varepsilon\ll\varepsilon$,
l'ensemble $\Gamma^\varepsilon_0$ est suffisamment grand pour
qu'\`a la limite on ait $\zeta_\alpha(0)=0$ pour $\alpha\in
\{2,3\}$. Quand $\varepsilon\ll r^\varepsilon\ll
\varepsilon^{1/3}$, on a en outre $\zeta_1(0)=0$. Finalement,
quand $\varepsilon^{1/3}\ll r^\varepsilon$, l'ensemble
$\Gamma^\varepsilon_0$ est si grand que toutes les conditions de
Dirichlet possibles sont satisfaites \`a la limite, c'est \`a dire
$\zeta_\alpha(0)=\zeta_1(0)= \frac{d\zeta_\alpha}{dy_1}(0)=c(0)=0$
(voir la version anglaise pour les d\'efinitions des espaces
fonctionnels qui interviennent dans le probl\`eme
limite~(\ref{18dicpbmalim}) 
et pour les d\'efinitions de ces fonctions).
\par Sauf quand $r^\varepsilon$ est de taille critique
(c'est \`a dire quand $r^\varepsilon\approx \varepsilon^\lambda$
avec $\lambda=3$, $1$ ou $1/3$), la forme bilin\'eaire ${\cal B}$
qui intervient dans (\ref{18dicpbmalim}) et la fonction
$P^\varepsilon$ qui intervient dans (\ref{correctore(Ueps)}) sont
nulles. Mais pour ces trois tailles critiques, ${\cal B}$ est une
forme bilin\'eaire coercive (en $\zeta_\alpha (0)$ pour
$\lambda=3$, en $\zeta_1(0)$ pour $\lambda=1$, et en $c(0)$ et
$\frac{d\zeta_\alpha}{d y_1}(0)$ pour $\lambda=1/3$, voir
(\ref{coerR}) de la version anglaise), de sorte que pour ces
tailles critiques $r^\varepsilon\approx\varepsilon^\lambda$,
${\cal B}$ est une 
p\'enalisation, par rapport au cas
$r^\varepsilon\ll \varepsilon^\lambda$,
des nouvelles conditions de
Dirichlet qui appara\^{i}tront pour
$r^\varepsilon \gg\varepsilon^\lambda$, ce qui r\'etablit une
certaine continuit\'e dans la transition.\par La pr\'esente Note
est la g\'en\'eralisation au cas de l'\'elasticit\'e lin\'eaire de
l'\'etude men\'ee dans \cite{CLMnc}, \cite{Casadobis}
pour le probl\`eme de diffusion analogue
quand le param\`etre
$t^\varepsilon$ qui intervient dans ces travaux est
$t^\varepsilon=0$. Les d\'emonstrations d\'etaill\'ees seront
donn\'ees dans \cite{Casado3}.

\setcounter{section}{0}
\selectlanguage{english}
\begin{figure}[t]
\begin{center}
$$\epsfig{file= 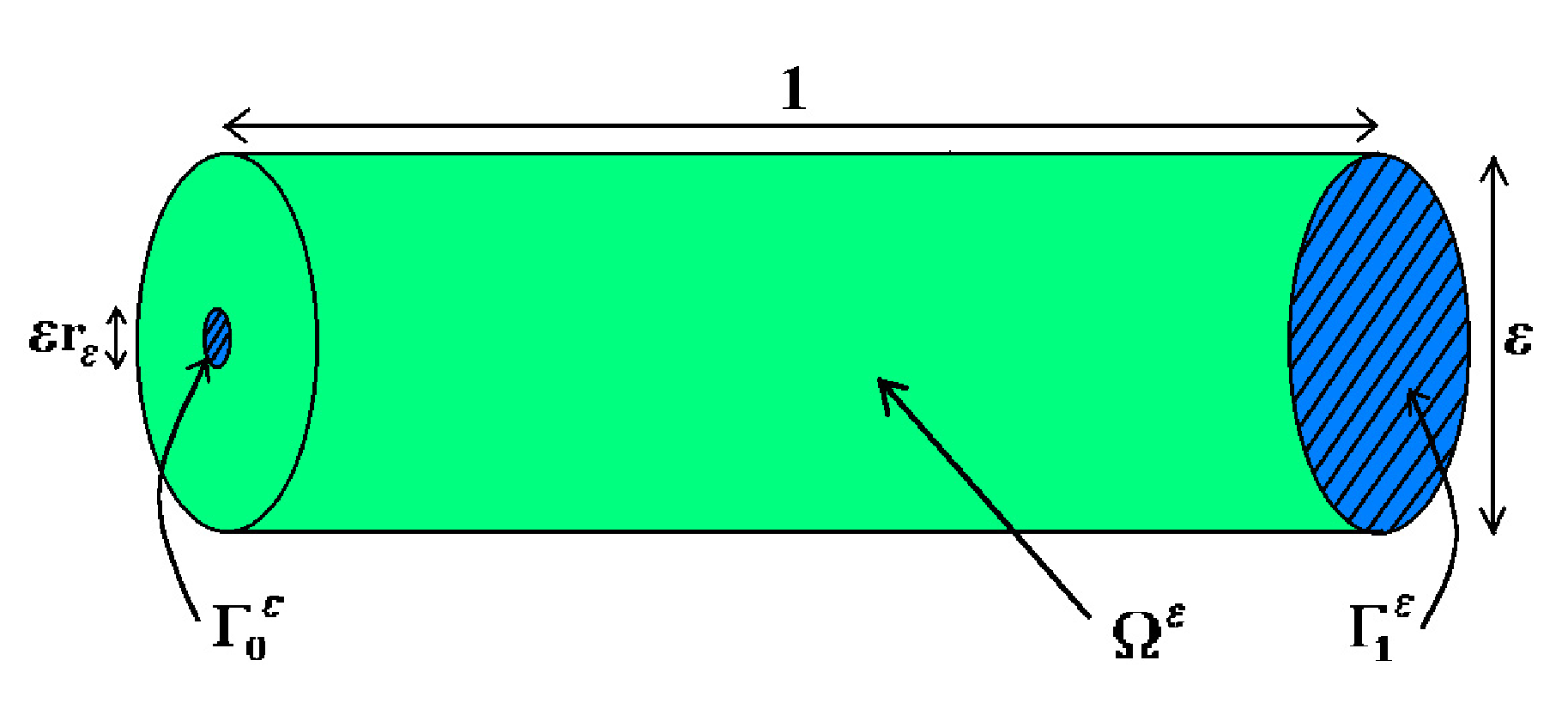,height=3.cm, width=7.5cm}$$
\vskip -0.8cm {\bf Figure 1.}
\end{center}
\vskip -0.1cm
\end{figure}
\section{Position of the problem and notation}In this Note we study
the asymptotic behavior of the solution of an anisotropic,
heterogeneous, linearized elasticity problem posed in a thin
cylinder~$\Omega^\varepsilon$ whose diameter $\varepsilon$ tends
to zero and whose axis is the first axis of coordinates~$Ox_1$. On
one of its extremities ($x_1=1$) the cylinder is fixed on its
whole basis~$\Gamma^\varepsilon_1$ whereas on the second one
($x_1=0$) it is fixed only on a small part~$\Gamma^\varepsilon_0$
of it, of diameter $\varepsilon r^\varepsilon$ much smaller than
$\varepsilon$. The Neumann boundary condition is assumed on the
remainder of the boundary of~$\Omega^\varepsilon$. Mathematically
the problem can be formulated as follows.
\par For $\varepsilon>0$, we consider $r^\varepsilon$ a positive parameter which tends to
zero with $\varepsilon$. Let $S_0$ and $S$ be two bounded smooth
domains in $\RR^2$, with $0\in S$. We define $\Omega=(0,1)\times
S$, $\Omega^\varepsilon=(0,1)\times \varepsilon S$ and
$\Gamma^\varepsilon=\Gamma^\varepsilon_0\cup
\Gamma^\varepsilon_1$, where
$\Gamma^\varepsilon_0=\{0\}\times\varepsilon r^\varepsilon S_0$,
$\Gamma^\varepsilon_1=\{1\}\times \varepsilon S$. Observe that the
size of $\Gamma^\varepsilon_0$ is much smaller than the size of
the basis $\{0\}\times \varepsilon S$ since $r^\varepsilon$
tends to zero. The thin cylinder $\Omega^\varepsilon$ is
represented in Figure~1 in the case where $S_0$ and $S$ are both
balls of $\RR^2$.
\par The elements of $\RR^3$ are decomposed as $x=(x_1,x')$,
$x_1\in\RR$, $x'=(x_2,x_3)\in\RR^2$. We denote by
$\{\mathbf{e}^1,\mathbf{e}^2,\mathbf{e}^3\}$ the canonical basis
of $\RR^3$ and by ${\cal L}(\RR^{3\times 3}_s)$ the space of
linear maps of $\RR^{3\times 3}_s$ into itself (or in other terms
of fourth order tensors), where $\RR^{3\times 3}_s$ is the space
of the $3\times 3$ symmetric matrices. We adopt Einstein's
convention of repeated indices. Greek indices ($\alpha$ and
$\beta$) take only the values $2$ and $3$, while latin indices
($i$ and $j$) take the values $1$, $2$ and $3$.
\par We consider $A\in C^0(\overline{\Omega};{\cal L}(\RR^{3\times 3}_s))$
such that there exists $m>0$ with\vskip -0.3cm $$A(y)\xi\xi \ge
m|\xi|^2,\quad\forall\xi\in\RR^{3\times 3}_s ,\quad\forall y\in
\overline{\Omega},$$\vskip -0.1cm\noindent and we define
$A^\varepsilon\in C^0(\Omega^\varepsilon;{\cal L}(\RR^{3\times
3}_s))$ by\vskip -0.2cm
\begin{equation}\label{defAeps}\displaystyle
A^\varepsilon(x)=A(x_1,\frac{x'}{\varepsilon}),\quad\forall x\in
\Omega^\varepsilon.\end{equation}\vskip -0.1cm
\noindent We also consider $f\in L^2(\Omega)^3$ and $h\in
L^2(\Omega)_s^{3\times 3}$, and we define $F^\varepsilon\in
L^2(\Omega^\varepsilon)^3$ and $H^\varepsilon\in
L^2(\Omega^\varepsilon)_s^{3\times 3}$ by\vskip -0.2cm
\begin{equation}\label{defFepsHeps}\displaystyle
F^\varepsilon(x)=f_1(x_1,\frac{x'}{\varepsilon})\mathbf{e}^1
+\varepsilon
f_\alpha(x_1,\frac{x'}{\varepsilon})\mathbf{e}^\alpha ,
\quad
H^\varepsilon(x)=h(x_1,\frac{x'}{\varepsilon}), \quad \mbox{a.e.
}x\in\Omega^\varepsilon.\end{equation}\vskip -0.1cm
\par In the thin domain $\Omega^\varepsilon$ we consider the
elasticity problem
\begin{equation}\label{probdehom}
\left\{\begin{array}{l} \displaystyle U^\varepsilon\in
H^1_{\Gamma^\varepsilon}(\Omega^\varepsilon)^3, \\
\noalign{\smallskip} \displaystyle \int_{\Omega^\varepsilon}
A^\varepsilon e(U^\varepsilon):e(\overline{U}^\varepsilon) dx =
\int_{\Omega^\varepsilon}  F^\varepsilon \overline{U}^\varepsilon
dx + \int_{\Omega^\varepsilon}
H^\varepsilon:e(\overline{U}^\varepsilon) dx,\quad\forall
\overline{U}^\varepsilon\in
H^1_{\Gamma^\varepsilon}(\Omega^\varepsilon)^3,\end{array}
\right.\end{equation}
where
$$\displaystyle
H^1_{\Gamma^\varepsilon}(\Omega^\varepsilon)=\left\{U\in
H^1(\Omega^\varepsilon)\,:\,U=0\mbox{ on
}\Gamma^\varepsilon\right\};$$
observe that in the above formulation, as well as in the remainder
of the present Note, complex numbers never appear, and that
$\overline{U}^\varepsilon$ (and later $\overline{u}$,
$\overline{v}$, $\overline{w}$) denotes the test function
associated to the solution $U^\varepsilon$, and not its complex
conjugate. Observe also that the solution $U^\varepsilon$
of~(\ref{probdehom}) satisfies a non homogeneous Neumann boundary
condition on the part $\partial\Omega^\varepsilon \setminus
\Gamma^\varepsilon$ where it is not fixed, since integrating by
parts $\int_{\Omega^\varepsilon}H^\varepsilon:e(U^\varepsilon)dx$
(when $h$ and therefore $H^\varepsilon$ is sufficiently smooth)
produces both body forces and surface forces. 
Similarly to the
body forces $F^\varepsilon$ we could have introduced explicit
surface forces $G^\varepsilon$ on
$\partial\Omega^\varepsilon\setminus \Gamma^\varepsilon$, but we
have preferred not to include them for the sake of simplicity.
\par It is well known that problem~(\ref{probdehom}) has an unique
solution (see, e.g., \cite{Ciarlet}). The aim of the present Note,
which announces our paper~\cite{Casado3}, is to describe the
asymptotic behavior of the solution $U^\varepsilon$ and to give a
corrector result for $e(U^\varepsilon)$ as $\varepsilon$ tends to
zero. The result depends on the behavior of $r^\varepsilon$ and
exhibits 3 critical sizes, 
namely $\varepsilon^3$, $\varepsilon$,
and $\varepsilon^{1/3}$, so that there are 7 different regimes:
$r^\varepsilon \ll \varepsilon^3$, $r^\varepsilon \approx
\varepsilon^3$, $\varepsilon^3 \ll r^\varepsilon \ll \varepsilon$,
$r^\varepsilon \approx \varepsilon$, $\varepsilon\ll r^\varepsilon
\ll \varepsilon^{1/3}$, 
$r^\varepsilon \approx \varepsilon^{1/3}$,
and $\varepsilon^{1/3} \ll r^\varepsilon \le C$, where
$r^\varepsilon \ll \varepsilon^\lambda$ stands for
$r^\varepsilon/\varepsilon^\lambda\to 0$ (and equivalently
$\varepsilon^\lambda \ll r^\varepsilon$ for
$r^\varepsilon/\varepsilon^\lambda\to +\infty$), while
$r^\varepsilon \approx \varepsilon^\lambda$ stands for
$r^\varepsilon/\varepsilon^\lambda\to \rho$, 
for some $\rho$ with
$0<\rho<+\infty$.
\par To express the results
and to make the proofs, we will use two changes of variables. The
first change of variables is given by
\par
\centerline{$\displaystyle y=y^\varepsilon(x)\quad\mbox{ with
}\quad\displaystyle y_1=x_1,\,\,
y'=\frac{x'}{\varepsilon},$}\vskip 0.1cm
\noindent which transforms the variable domain
$\Omega^\varepsilon$ into the fixed domain~$\Omega$. This is the
usual change of variables used to study equations in thin
cylinders (see, e.g., \cite{Cimetiere}, \cite{MuratSiliE},
\cite{MuratSilitoappear}, \cite{Trabucho}). When $r^\varepsilon=1$
and $S_0=S$ (but the same proof works for $r^\varepsilon=C$
independent of $\varepsilon$ such that $CS_0\subset S$), it was
used successfully 
in~ \cite{Cimetiere}, \cite{MuratSiliE}, \cite{MuratSilitoappear}
to pass to the limit in~(\ref{probdehom}). But when
$r^\varepsilon$ tends to zero with~$\varepsilon$, this first
change of variables does not provide the information we need about
the behavior of $U^\varepsilon$ in the part of
$\Omega^\varepsilon$ close to $\Gamma^\varepsilon_0$. Thus we
introduce a second change of variables given by
\par
\centerline{$\displaystyle z=z^\varepsilon(x)\quad\mbox{ with
}\quad\displaystyle z=\frac{x}{\varepsilon r^\varepsilon},$}\vskip
0.1cm
\noindent which transforms the variable domain
$\Omega^\varepsilon$ into a variable domain~$Z^\varepsilon$, the
limit of which is the half space $Z=(0,+\infty)\times\RR^2$.
Observe that the Dirichlet boundary condition is now imposed on
the fixed part~$\{0\}~\times~S_0$ 
of the boundary of~$Z^\varepsilon$. This
change of variables provides a suitable rescaling near~$x_1=0$. It
was used successfully in~\cite{CLMnc}, \cite{Casadobis} to study
the diffusion problem similar to~(\ref{probdehom}) in the geometry
that we consider in the present Note (and even in a more
complicated one, where $\Omega^\varepsilon$ is made of union of
two cylinders $\displaystyle \left\{(-t^\varepsilon,0)\times
\varepsilon r^\varepsilon S_0\right\}\cup \left\{(0,1)\times
\varepsilon S\right\}$ and where the Dirichlet condition is
imposed on $x_1=-t^\varepsilon$ and $x_1=1$; the geometry
considered in the present Note corresponds to $t^\varepsilon=0$; a
problem of conduction in a notched beam of the same type was
solved in~\cite{Casado} by using the same change of variables).
\par We denote by $D^{1,2}(Z)$ the Deny's space
$$D^{1,2}(Z)=\left\{p:p\in
L^6(Z),\,\,\nabla p\in L^2(Z)^3\right\}.$$
We will also use the functional spaces (already used
in~\cite{Gaudiello}, \cite{Gaudiello2}, \cite{MuratSiliE},
\cite{MuratSilitoappear})
$$\left\{\begin{array}{ll}\displaystyle BN_b(\Omega)=\Big\{u: & \displaystyle
\exists \zeta_\alpha\in
H^2(0,1),\,\,\zeta_\alpha(1)=\frac{d\zeta_\alpha}{dy_1}(1)=0,\,\,
u_\alpha(y)=\zeta_\alpha(y_1), \,\,\forall\alpha\in\{2,3\}, \\
\noalign{\smallskip} & \displaystyle \exists \zeta_1 \in
H^1(0,1),\,\,\zeta_1(1)=0,\,\,u_1(y)=\zeta_1(y)-\frac{d\zeta_\alpha}{dy_1}(y_1)y_\alpha
\Big\},\end{array}\right.$$
$$\left\{\begin{array}{ll}\displaystyle R_b(\Omega)=\Big\{v: & \displaystyle v_1\in L^2(0,1;H^1(S)),
\,\,\int_{S}v_1(y_1,y')dy'=0\mbox{ a.e. }y_1\in(0,1), \\
\noalign{\smallskip}\displaystyle & \exists c\in
H^1(0,1),\,\,c(1)=0,\,\,
v_2(y)=c(y_1)y_3,\,\,v_3(y)=-c(y_1)y_2\Big\},\end{array}\right.$$
$$\left\{\begin{array}{ll}\displaystyle RD_2^\perp(\Omega)=
\Big\{w: & \displaystyle w_1=0,\,\,w_\alpha\in L^2(0,1;H^1(S)),\,
\displaystyle \int_S w_\alpha(y_1,y')dy'=0,
\,\,\forall \alpha\in\{2,3\},
\\
\noalign{\smallskip} & \displaystyle \int_S\left(y_3w_2(y_1,y')-
y_2w_3(y_1,y')\right) dy'\hskip -0.1cm= 0 \,\,\mbox{ a.e. }y_1\in
(0,1) \Big\}.\end{array}\right.$$
For a given $(u,v,w)\in BN_b(\Omega)\times R_b(\Omega) \times
RD^\perp_2(\Omega)$, we denote by $E(u,v,w)$ the second order
symmetric tensor with values in $\RR^{3\times 3}_s$ defined by
\par
\vskip 0.1cm
\centerline{$\displaystyle
E_{11}(u,v,w)=e_{11}(u),\,\,E_{1\beta}(u,v,w)=e_{1\beta}(v),\,\,E_{\alpha\beta}(u,v,w)=e_{\alpha\beta}(w),\,\,\forall
\alpha,\beta\in\{2,3\}.$}

\section{The result and some comments}

The asymptotic behavior of the solution of~(\ref{probdehom})
depends on the size of $r^\varepsilon$ with respect to
$\varepsilon$. Seven regimes appear in the following Theorem which
describes the asymptotic behavior of~$U^\varepsilon$ and provides
a corrector result for~$U^\varepsilon$ and~$e(U^\varepsilon)$.
\vskip 0.1cm
\par {\sc Theorem\, \rule[0.75mm]{1.75mm}{0.4pt}\, }{\it Let $U^\varepsilon$ be the solution
of~(\ref{probdehom}). There exist a closed linear subspace ${\cal E}$ 
of $BN_b(\Omega)\times R_b(\Omega)\times RD^\perp_2(\Omega)$,
a function $P^\varepsilon\in L^2(\Omega^\varepsilon)_s^{3\times
3}$, and a bilinear continuous nonnegative form ${\cal B}$ on
${\cal E}\times {\cal E}$ such that, 
defining $(u,v,w)$ as the solution
of the variational problem
\begin{equation}\label{18dicpbmalim}\left\{\begin{array}{l} 
\displaystyle (u,v,w)\in{\cal E}, \\
\noalign{\smallskip} \displaystyle \int_\Omega
AE(u,v,w):E(\overline{u},\overline{v},\overline{w})dy +
\langle\,
{\cal B}(u,v,w),(\overline{u},\overline{v},\overline{w})\,\rangle =
\\
\noalign{\smallskip} \displaystyle = \int_\Omega f\overline{u}dy +
\int_\Omega
h:E(\overline{u},\overline{v},\overline{w})dy,\quad\forall
(\overline{u},\overline{v},\overline{w})\in{\cal E},
\end{array}\right.\end{equation}
then, when $\varepsilon$ tends to zero, we have
\begin{eqnarray}\displaystyle\frac{1}{|\Omega^\varepsilon|}
\int_{\Omega^\varepsilon}\left(\,|U^\varepsilon_1(x)-
u_1(x_1,\frac{x'}{\varepsilon})|^2+\sum_{\alpha=2}^3|\varepsilon
U^\varepsilon_\alpha(x)-u_\alpha(x_1) |^2\,\right)dx
\longrightarrow 0, \label{correctorUeps} \\ \noalign{\smallskip}
\displaystyle\frac{1}{|\Omega^\varepsilon|}\int_{\Omega^\varepsilon}|e(U^\varepsilon)(x)
- E(u,v,w)(x_1,\frac{x'}{\varepsilon}) -
P^\varepsilon(\frac{x}{\varepsilon r^\varepsilon})|^2dx
\longrightarrow 0.\label{correctore(Ueps)}\end{eqnarray}
\par The definitions of ${\cal E}$, $P^\varepsilon$ and ${\cal B}$
do not depend on the forces $f$ and $h$ which define
$F^\varepsilon$ and $H^\varepsilon$, but only on the set~$S_0$, on
the fourth order tensor~$A$, and on the behavior of
$r^\varepsilon$ when $\varepsilon$ tends to zero, and more
specifically of its behavior with respect to $\varepsilon^3$,
$\varepsilon$, and $\varepsilon^{1/3}$, such that there are 7
regimes, which are described now.
\begin{itemize}
\item If $r^\varepsilon \ll \varepsilon^3$, then
${\cal E}=BN_b(\Omega)\times R_b(\Omega)\times
RD^\perp_2(\Omega),$
$P^\varepsilon=0$, ${\cal B}=0$.
\item If $r^\varepsilon \approx
\varepsilon^3$ with $r^\varepsilon/\varepsilon^3\to \rho$, then
${\cal E}=BN_b(\Omega)\times R_b(\Omega)\times
RD^\perp_2(\Omega)$, and defining $\varphi^\alpha$,
$\alpha\in\{2,3\}$, as the solution of
$$\left\{\begin{array}{l}
\displaystyle \varphi^\alpha\in D^{1,2}(Z)^3,\,\,
\varphi^\alpha= \mathbf{e}^\alpha\mbox{ on }\{0\}\times S_0, \\
\noalign{\smallskip} \displaystyle \int_Z
A(0)e(\varphi^\alpha):e(\overline{\varphi})dz=0,\quad\forall
\overline{\varphi}\in D^{1,2}(Z)^3,\,\,\overline{\varphi}=0\mbox{
on }\{0\}\times S_0,
\end{array}\right.$$
then one has
$$P^\varepsilon(z)=-\frac{1}{\varepsilon^2 r^\varepsilon}\zeta_\alpha(0)
\,e(\varphi^\alpha)(z),\,\,\mbox{ a.e. }z\in Z,$$
$$\langle\,{\cal B}(u,v,w),(\overline{u},\overline{v},\overline{w})\,\rangle=\rho
\int_Z
A(0)\left(\zeta_\alpha(0)\,e(\varphi^\alpha)\right):
\left(\overline{\zeta}_\beta(0)\,e(\varphi^\beta)\right)dz,
\,\,\forall
(u,v,w),(\overline{u},\overline{v},\overline{w})\in{\cal E}.$$
\item If $\varepsilon^3 \ll r^\varepsilon \ll \varepsilon$, then
$$\displaystyle {\cal E}=\left\{(u,v,w)\in BN_b(\Omega)\times
R_b(\Omega)\times
RD^\perp_2(\Omega):\zeta_\alpha(0)=0,\,\,\forall\alpha\in\{2,3\}\right\},
\,\,P^\varepsilon=0,\,\,{\cal B}=0.$$
\item If $r^\varepsilon \approx \varepsilon$ with $r^\varepsilon/\varepsilon \to \rho$, then
$$\displaystyle {\cal E}=\left\{(u,v,w)\in BN_b(\Omega)\times
R_b(\Omega)\times
RD^\perp_2(\Omega):\zeta_\alpha(0)=0,\,\,\forall\alpha\in\{2,3\}\right\},$$
and defining $\varphi^1$ as the solution of
$$\left\{\begin{array}{l}
\displaystyle \varphi^1\in D^{1,2}(Z)^3,\,\,
\varphi^1= \mathbf{e}^1\mbox{ on }\{0\}\times S_0, \\
\noalign{\smallskip} \displaystyle \int_Z
A(0)e(\varphi^1):e(\overline{\varphi})dz=0,\quad\forall
\overline{\varphi}\in D^{1,2}(Z)^3,\,\,\overline{\varphi}=0\mbox{
on }\{0\}\times S_0,
\end{array}\right.$$
and then setting
$\displaystyle\hat{\varphi}^1=\varphi^1+a_\alpha\,\varphi^\alpha$,
where $(a_2,a_3)$ is defined by
$$\left\{\begin{array}{l}
\displaystyle (a_2,a_3)\in\RR^2, \\
\noalign{\smallskip} \displaystyle \int_Z A(0)\left(e(\varphi^1)+
a_\alpha\, e(\varphi^\alpha)\right) :\left(\overline{a}_\beta\,
e(\varphi^\beta)\right)dz=0,\quad\forall
(\overline{a}_2,\overline{a}_3)\in\RR^2,\end{array}\right.$$
then one has
$$\displaystyle P^\varepsilon(z)=-\frac{1}{\varepsilon r^\varepsilon}\,\zeta_1(0)\,
e(\hat{\varphi}^1)(z), \mbox{ a.e. }z\in Z,$$
$$\langle\,{\cal B}(u,v,w),(\overline{u},\overline{v},\overline{w})\,\rangle=\rho
\int_Z A(0)\left(\zeta_1(0)\, e(\hat{\varphi}^1)\right) :
\left(\overline{\zeta}_1(0)\,
e(\hat{\varphi}^1)\right)dz,\,\forall
(u,v,w),(\overline{u},\overline{v},\overline{w})\in{\cal E}.$$
\item If $\varepsilon\ll r^\varepsilon
\ll \varepsilon^{1/3}$, then
$$\displaystyle {\cal E}=\left\{(u,v,w) \hskip -0.05cm \in \hskip -0.05cm BN_b(\Omega)
\hskip -0.05cm \times \hskip -0.05cm R_b(\Omega) \hskip -0.05cm
\times \hskip -0.05cm RD^\perp_2(\Omega):\zeta_\alpha(0) \hskip
-0.1cm = \hskip -0.1cm \zeta_1(0) \hskip -0.1cm = \hskip -0.1cm
0,\,\forall \alpha\hskip -0.1cm \in\hskip -0.1cm
\{2,3\}\right\},\,P^\varepsilon \hskip -0.05cm = \hskip -0.05cm
0,\,{\cal B} \hskip -0.05cm = \hskip -0.05cm 0.$$
\item If $r^\varepsilon \approx \varepsilon^{1/3}$ with
$(r^\varepsilon)^3/\varepsilon\to \rho$, then
$$\displaystyle {\cal E}=\left\{(u,v,w)\in BN_b(\Omega)\times
R_b(\Omega)\times
RD^\perp_2(\Omega):\zeta_\alpha(0)=\zeta_1(0)=0,\,\,\forall
\alpha\in\{2,3\}\right\},$$
and defining $\psi^1$ as the solution of
$$\left\{\begin{array}{l} \displaystyle
\psi^1 \in D^{1,2}(Z)^3,\,\,
\psi^1=z_3\mathbf{e}^2 - z_2\mathbf{e}^3\mbox{ on }\{0\}\times S_0, \\
\noalign{\smallskip} \displaystyle \int_Z
A(0)e(\psi^1):e(\overline{\psi})dz=0,\quad\forall
\overline{\psi}\in D^{1,2}(Z)^3,\,\,\overline{\psi}=0\mbox{ on
}\{0\}\times S_0,
\end{array}\right.$$
and $\psi^\alpha$, $\alpha\in\{2,3\}$, as the solution of
$$\left\{\begin{array}{l} \displaystyle
\psi^\alpha \in D^{1,2}(Z)^3,\,\, \psi^\alpha=
z_1\mathbf{e}^\alpha - z_\alpha\mathbf{e}^1
\mbox{ on }\{0\}\times S_0, \\
\noalign{\smallskip} \displaystyle \int_Z
A(0)e(\psi^\alpha):e(\overline{\psi})dz=0,\quad\forall
\overline{\psi}\in D^{1,2}(Z)^3,\,\, \overline{\psi}=0\mbox{ on
}\{0\}\times S_0,
\end{array}\right.$$
and then setting
$\displaystyle\hat{\psi}^i=\psi^i+b^i_k\,\varphi^k$, where
$(b^i_1,b^i_2,b^i_3)$, $i\in\{1,2,3\}$, is defined by
$$\left\{\begin{array}{l}
\displaystyle (b^i_1,b^i_2,b^i_3)\in\RR^3, \\ \noalign{\smallskip}
\displaystyle \int_Z A(0)\left(e(\psi^i)+b^i_k\,
e(\varphi^k)\right):\left(\overline{b}^i_l\,
e(\varphi^l)\right)dz=0,\quad
\forall(\overline{b}^i_1,\overline{b}^i_2,\overline{b}^i_3)\in\RR^3,
\end{array}\right.$$
then one has
$$P^\varepsilon(z)=-\frac{1}{\varepsilon}\left(c(0)\,e(\hat{\psi}^1)
(z) + \frac{d \zeta_\alpha}{d y_1}(0)\,e(\hat{\psi}^\alpha)(z)
\right),\mbox{ a.e. }z\in Z,$$
$$\left\{\begin{array}{l}\displaystyle\langle\,{\cal B}(u,v,w),
(\overline{u},\overline{v},\overline{w})\,\rangle=\rho \int_Z
A(0)\left(c(0)\, e(\hat{\psi}^1)+\frac{d \zeta_\alpha}{d
y_1}(0)\,e(\hat{\psi}^\alpha)\right):
\\
\noalign{\smallskip} \displaystyle \hskip 2.cm
:\left(\overline{c}(0)\, e(\hat{\psi}^1)+\frac{d
\overline{\zeta}_\alpha}{d y_1}(0)\,
e(\hat{\psi}^\alpha)\right)dz,\quad \forall
(u,v,w),(\overline{u},\overline{v},\overline{w})\in{\cal E}.
\end{array}\right.$$
\item If $\varepsilon^{1/3} \ll r^\varepsilon \le C$, then
$$\displaystyle \begin{array}{c}\displaystyle {\cal E}=\Big\{(u,v,w)\hskip
-0.05cm\in\hskip -0.05cm BN_b(\Omega)\hskip -0.05cm\times\hskip
-0.05cm R_b(\Omega)\hskip -0.05cm\times\hskip -0.05cm
RD^\perp_2(\Omega): \zeta_\alpha(0)\hskip -0.1cm=\zeta_1(0)\hskip
-0.1cm=\frac{d\zeta_\alpha}{dy_1}(0)\hskip -0.1cm =c(0)\hskip
-0.1cm=0,\,\forall\alpha\hskip -0.05cm\in\hskip
-0.05cm\{2,3\}\Big\},\end{array}$$
$P^\varepsilon=0$, ${\cal B}=0$.
\end{itemize}}
\vskip 0.1cm
Let us make some comments about the statement of this Theorem.
\par Assertion~(\ref{correctore(Ueps)}) is a corrector result, since
using $\overline{U}^\varepsilon=U^\varepsilon$ as test function
in~(\ref{probdehom}) and Korn's inequality, one can prove that
\par\vskip 0.2cm
\centerline{$\displaystyle\frac{1}{|\Omega^\varepsilon|}\int_{\Omega^\varepsilon}|e(U^\varepsilon)|^2dx
\le C\left[\frac{1}{|\Omega|}\int_\Omega |f|^2dx +
\frac{1}{|\Omega|}\int_\Omega |h|^2dx\right]$,}\vskip 0.2cm
\noindent and since one can also prove that
\par\vskip 0.2cm
\centerline{$\displaystyle
\frac{1}{|\Omega^\varepsilon|}\int_{\Omega^\varepsilon}|E(u,v,w)(x_1,\frac{x'}{\varepsilon})|^2dx
+
\frac{1}{|\Omega^\varepsilon|}\int_{\Omega^\varepsilon}|P^\varepsilon(\frac{x}{\varepsilon
r^\varepsilon})|^2dx \approx 1.$}\vskip 0.2cm
\par If one examines the definition of ${\cal E}$, one realizes that
the number of Dirichlet boundary conditions imposed in the
definition of ${\cal E}$  increases with the size
of~$r^\varepsilon$. Indeed, in view of the definitions of
$BN_b(\Omega)$, $R_b(\Omega)$ and $RD^\perp_2(\Omega)$, the sole
functions which have a trace for $x_1=0$ are $\zeta_\alpha$,
$\zeta_1$, $\frac{d\zeta_\alpha}{dy_1}$ and $c$, where
$\alpha\in\{2,3\}$ (the other functions, namely $v_1$ and
$w_\alpha$, have no trace for $x_1=0$). When $r^\varepsilon\ll
\varepsilon^3$, the set
$\Gamma^\varepsilon_0=\{0\}\times\varepsilon r^\varepsilon S_0$ is
too small and the homogeneous Dirichlet boundary condition imposed
for $x_1=0$ to the solution $U^\varepsilon$ of~(\ref{probdehom})
completely disappears at the limit. When $\varepsilon^3\ll
r^\varepsilon\ll \varepsilon$, the set $\Gamma^\varepsilon_0$ is
sufficiently large to impose at the limit that $\zeta_\alpha(0)=0$
for $\alpha\in\{2,3\}$. When $\varepsilon\ll r^\varepsilon\ll
\varepsilon^{1/3}$, one further has $\zeta_1(0)=0$. Finally when
$\varepsilon^{1/3}\ll r^\varepsilon$, the set
$\Gamma^\varepsilon_0$ is so large that all the possible Dirichlet
boundary conditions are imposed at $x_1=0$, namely
$\zeta_\alpha(0)=\zeta_1(0)=\frac{d\zeta_\alpha}{dy_1}(0)=c(0)=0$.
\par Except in the three regimes where the size of $r^\varepsilon$
is critical (i.e. when $r^\varepsilon\approx \varepsilon^\lambda$,
with $\lambda=3$, $1$, or $1/3$), one always has $P^\varepsilon=0$
and ${\cal B}=0$. For these three critical sizes, one can show
that ${\cal B}$ is a coercive bilinear form, in the sense that
there exists some $n>0$ such that\vskip -0.35cm
\begin{equation}\label{coerR}\left\{\begin{array}{ll} \displaystyle
\langle\,{\cal B}(u,v,w),(u,v,w)\,\rangle\ge n\rho
\sum_{\alpha=2,3}|\zeta_\alpha(0)|^2, & \displaystyle \mbox{when
$r^\varepsilon\approx \varepsilon^3$,}
\\ \noalign{\smallskip} \displaystyle \langle\,{\cal B}(u,v,w),(u,v,w)\,\rangle\ge n\rho
|\zeta_1(0)|^2, & \displaystyle \mbox{when $r^\varepsilon\approx
\varepsilon$,}
\\
\noalign{\smallskip} \displaystyle \langle\,
{\cal B}(u,v,w),(u,v,w)\,\rangle\,\ge n\rho \left( |c(0)|^2 +
\sum_{\alpha=2,3}\left|\frac{d\zeta_\alpha}{dy_1}(0)\right|^2
\right), & \displaystyle \mbox{when $r^\varepsilon\approx
\varepsilon^{1/3}$.} \end{array}\right.\end{equation}\vskip -0.2cm
\noindent This implies that for every critical size $r^\varepsilon
\approx \varepsilon^\lambda$, the new Dirichlet boundary
conditions which appear for $r^\varepsilon \gg
\varepsilon^\lambda$ (with respect to those imposed for
$r^\varepsilon \ll \varepsilon^\lambda$) are penalized by the
value of~$\rho$. This introduces some type of continuity in the
transition of the Dirichlet condition between the two regimes
which are separated by a critical size~$\varepsilon^\lambda$. For
these three critical sizes, the functions $\varphi^\alpha$,
$\varphi^1$, $\hat{\varphi}^1$, $\psi^i$, and $\hat{\psi}^i$ are
in some sense generalized capacitary potentials 
of $\{0\}\times
S_0$ in $Z$, and the bilinear form ${\cal B}$ is in some sense an
asymptotic trace of some type of capacity of
$\Gamma^\varepsilon_0$ in $\Omega^\varepsilon$ for the weighted
energy
$\frac{1}{|\Omega^\varepsilon|}\int_{\Omega^\varepsilon}A(x)e(\varphi):e(\varphi)dx$.
\par The present work is the natural generalization to the
elastic case of~\cite{CLMnc}, \cite{Casadobis}, where diffusion
problems were posed in the union of two cylinders
$\displaystyle\left\{(-t^\varepsilon,0)\times \varepsilon
r^\varepsilon S_0\right\}\cup \left\{(0,1)\times \varepsilon
S\right\}$, with both $t^\varepsilon$ and $r^\varepsilon$ tending
to zero (the present geometry corresponds to $t^\varepsilon=0$).
When $t^\varepsilon=0$, the diffusion problem was in
comparison more simple, since only one critical 
size, namely $r^\varepsilon\approx \varepsilon$, 
appeared in the analysis,
separating a pure Neumann boundary condition [corresponding to the
analogue of $u$ satisfying $u\in H^1(0,1)$, $u(1)=0$] for
$r^\varepsilon\ll\varepsilon$ and a pure Dirichlet boundary
condition [corresponding to the analogue of $u$ satisfying $u\in
H^1(0,1)$, $u(0)=u(1)=0$] for $r^\varepsilon \gg \varepsilon$.
These works were related to~\cite{Casado}, where a notched beam
for diffusion problems was considered. The present work is also
related to~\cite{Gaudiello}, \cite{Gaudiello2}, where a
multidomain made of an elastic vertical beam of length~1 and of
radius~$r^\varepsilon$ and of an horizontal plate of radius~1 and
of height~$\varepsilon$ was considered.
\par The detailed proofs of the results of the present Note will be given in a
forthcoming paper~\cite{Casado3}.
\Acknowledgements{This work has been partly supported by the
project BFM2002-00672 of the DGI of Spain and by the HMS2000
Training and Research Network. It was begun during a visit of the
third author to the Departamento de Ecuaciones Diferenciales y
An\'alisis Num\'erico de la Universidad de Sevilla.}

%

\begin{thebibliography}{99}
\selectlanguage{english}
\bibitem{Casado}  {Casado-D\'{\i}az J., Murat F.,} The diffusion
equation in a notched beam, to appear.
%
\bibitem{CLMnc} {Casado-D\'{\i}az J., Luna-Laynez M., Murat
F.,} Asymptotic behavior of diffusion problems in a domain made of
two cylinders of different diameters and lengths, C. R. Acad. Sci.
Paris, Ser. I, 338 (2004), 133-138.
%
\bibitem{Casadobis} {Casado-D\'{\i}az J., Luna-Laynez M., Murat
F.,} Diffusion problems in a domain made of two thin cylinders of
different diameters and lengths, to appear.
%
\bibitem{Casado3} {Casado-D\'{\i}az J., Luna-Laynez M., Murat
F.,} Elasticity problems in a beam fixed on only a small part of
one of its extremities, to appear.
%
\bibitem{Ciarlet} {Ciarlet P.G.,} Mathematical elasticity, Vol. I:
Three-dimensional elasticity, North-Holland, 1988.
%
\bibitem{Cimetiere}
{Cimeti\`ere A., Geymonat G., Le Dret H., Raoult A., Tutek Z.,}
Asymptotic theory and analysis for displacements and stress
distribution in nonlinear straight slender rods, J. of Elasticity,
19 (1988), 111-161.
%
\bibitem{Gaudiello} {Gaudiello A., Monneau R., Mossino J., Murat F., Sili A.,} On the junction of
elastic plates and beams, C. R. Acad. Sci. Paris, Ser. I, 335
(2002), 717-722.
\bibitem{Gaudiello2} {Gaudiello A., Monneau R., Mossino J., Murat F., Sili
A.,} Junction of elastic plates and beams, to appear.
%
\bibitem{MuratSiliE}  {Murat F., Sili A.,}
Comportement asymptotique des solutions du syst\`{e}me de
l'\'elasticit\'e lin\'earis\'ee anisotrope h\'et\'erog\`ene dans
des cylindres minces, C. R. Acad. Sci. Paris, Ser. I, 328 (1999),
179-184.
%
\bibitem{MuratSilitoappear}  {Murat F., Sili A.,} Anisotropic, heterogeneous, linearized
elasticity in thin cylinders, to appear.
%
\bibitem{Trabucho} {Trabucho L., Via\~no J.M.,} Mathematical modelling
of rods. Handbook of Numerical Analysis, Vol. IV, North-Holland,
1996.
%
\end{thebibliography}
\end{document}